\documentclass[reqno]{amsart}

\usepackage[utf8]{inputenc}
\usepackage{hyperref}

\usepackage{changepage}
\usepackage[pdftex]{graphicx}
\usepackage{color}
\usepackage{rotating} 

\usepackage{amsmath}
\usepackage{amsfonts}
\usepackage{amssymb}
\usepackage{amsthm}
\usepackage{mathtools}

\usepackage{framed}
 \usepackage{paralist}
\usepackage[shortlabels]{enumitem}

\usepackage{thmtools}
\usepackage{thm-restate}

\usepackage{hyperref}

\usepackage{cleveref}

\newtheorem{thm}{Theorem}

\newtheorem{Def}{Definition}

\newcommand\ringZ{\mathbb{Z}}       
   
\newcommand\ringZgeq{\mathbb{Z}_{\geq 0}}           
\newcommand\ringR{\mathbb{R}}  
  
\newcommand\ringRgeq{\mathbb{R}_{\geq 0}}              
\newcommand\ringX{\mathcal{X}}
\newcommand\ringvol{\textnormal{vol}}
\newcommand\ringE{\textnormal{E}_{P}}
\newcommand\ringG{\mathcal{G}}
\newcommand\ringw{correction volume }

\newcommand{\ringfcone}[2]{\textnormal{fcone}(#1,#2)}
\newcommand{\ringRS}{RS-$\mu$-value }
\newcommand{\ringRSs}{RS-$\mu$-values }
\newcommand{\ringRSsc}{RS-$\mu$-values}
\newcommand{\ringRSC}{RS-$\mu$-construction }
\newcommand{\ringRSCc}{RS-$\mu$-construction}
\newcommand{\ringStrip}[1]{\textnormal{Strip}(#1)}
\newcommand{\ringReg}{\tilde R}
\newcommand{\ringlin}{\ringlinsub}
\newcommand{\ringInt}{\textnormal{int}}
\newcommand{\ringconv}{\textnormal{conv}}
\newcommand{\ringlinsub}{\textnormal{lineal}}
\newcommand{\ringlat}{\textnormal{lat}}
\newcommand\ringsymm{\sigma_0}
\newcommand{\ringLat}{\ringlat}

\begin{document}

\title{Special cases and a dual view on the local formulas for Ehrhart coefficients from lattice tiles}

\author{M.~H.~Ring}

\address{Institute for Mathematics, 
University of Rostock,\\
18057 Rostock,
Germany\\
E-mail: maren.ring@uni-rostock.de\\
https://www.geometrie.uni-rostock.de/people/maren/}

\begin{abstract}
McMullen's formulas or local formulas for Ehrhart coefficients are functions on rational cones that determine the $i$-th coefficient of the Ehrhart polynomial as a weighted sum of the volumes of the i-dimensional faces of a polytope. 
This work focuses on the \ringRSC as given in \cite{RS17}. We give an explicit description of the construction from the dual point of view, i.e. given the cone of feasible directions instead of the normal cone as input value. We further show some properties of the construction in special cases, namely in case of symmetry and for the codimension one case.
\end{abstract}

\keywords{Ehrhart coefficients, local formula, lattice tiling.} 

\maketitle

%\bodymatter

\section{Introduction}

For a lattice polytope $P$ in a Euclidean space $V$ with lattice $\Lambda$, the \emph{Ehrhart polynomial} $\ringE$ counts the number of lattice points in the $t$-th dilate of $P$, 
\begin{equation*}
\ringE(t)=|\Lambda \cap t P|=e_dt^d+e_{d-1}t^{d-1}+\ldots + e_1t+e_0,
\end{equation*} 
for all $t\in \ringZgeq$ and with $d=\dim(P)$ (cf. Ehrhart 1962 \cite{Ehrhart}).

In 1983, McMullen \cite{mcmullen83} showed the existence of so-called \emph{McMullen's formulas} or \emph{local formulas} for Ehrhart polynomials. 
\begin{Def} \label{Def:LocalFormula}
A real valued function $\mu$ on rational cones in $V$ is called a \emph{McMullen's formula} or \emph{local formula} for Ehrhart coefficients, if for any lattice polytope $P$ with Ehrhart polynomial $\ringE (t)=e_d t^d+e_{d-1}t^{d-1}+\dots + e_1t+e_0$, we have
\begin{align*}
e_i=\sum_{\mathcal{F}_i} \mu(N_f) \ringvol(f), 
\end{align*}
for all $i\in\{0,\dots, d\}$ with $\mathcal{F}_i$ the set of all $i$-dimensional faces of $P$. 
\end{Def}
Here, $N_f$ is the \emph{(outer) normal vector} of the face $f$, i.e. the cone over the outer normal vectors of all facets meeting in $f$. It is a \emph{convex rational cone}, which means that it can be written as $\{a_1 v_1 + \ldots + a_k v_k \mid a_i \in \ringRgeq \}$ for some vectors $a_1, \ldots, a_k \in \Lambda$. The volume $\ringvol(f)$ denotes the \emph{relative volume} of $f$ with respect to the induced lattice in the affine span of $f$, as defined in Section \ref{Sec:Construction}. 

There have been several nice constructions of McMullen's formulas, for example the ones by Pommersheim and Thomas, see \cite{PT04}, and Berline and Vergne, see \cite{BV07}. This work will focus on the construction given in \cite{RS17} by the author together with Achill Sch\"urmann. Following the notation of Federico Castillo in his talk at Osaka University in 2018, we will refer to our construction as the \ringRSCc, in analogy to the BV-$\alpha$-construction by Berline and Vergne (cf. \cite{CL15}). 
Each construction has different advantages. This construction, for instance, does not require a prior decomposition into unimodular cones and can be described using elementary geometric means. 
In contrast to \cite{RS17}, it might sometimes be more convenient to consider $\mu$ not as a function on the normal cones of a polytope but from the dual point of view, given \emph{cones of feasible directions} as input values. Note that both cones, the normal and the cone of feasible directions of a face, hold the same local information about the face and can be recovered from one another. In Section \ref{Sec:Construction} we give a detailed construction of the \ringRSs given a full dimensional cone. 
To give more insight on how to compute the values, we then follow the construction given a polytope with cones of feasible directions and alongside compute the values in  an example. 
In Section \ref{Sec:Motivation} we focus on the connection between the \ringRSs and the number of lattice points, thus giving a motivation for the idea behind the construction and an outline of the proof. 
In Section \ref{Sec:SpecialCases} we focus on symmetry and show a new result about the codimension one case under central symmetry.

\section{Construction}\label{Sec:Construction}

Fix an ambient Euclidean space $V$ and a lattice $\Lambda$ of full rank.
Let $Q$ be a polyhedron and $f$ a face of $Q$. We define the \emph{cone of feasible directions} of $Q$ in the face $f$, \emph{fcone} of $Q$ in $f$ for short, as  
\begin{equation*}
\ringfcone Qf=\{x \in V \mid \exists \varepsilon>0 \colon s+\varepsilon x \in Q\}
\end{equation*} 
for any point $s$ in the relative interior $\ringInt(f)$ of $f$.
 Let further $\ringlinsub(Q)$ be the \emph{lineality space} of $Q$, i.e. the biggest linear subspace contained in $Q$, and denote $\ringlat(Q):= \Lambda\cap \ringlinsub(Q)$.

The construction we will give relies on a choice of certain lattice tiles, namely \emph{fundamental domains} as defined below. Different choices lead to different values and thus give an infinite family of constructions. 
\begin{Def}
For a polyhedron $Q\subseteq V$ with induced sublattice $L= \ringlinsub(Q)\cap \Lambda$ in the lineality space of $Q$, a \emph{fundamental domain} $T(Q)$ is a subset of $\ringlinsub(Q)$ with the following properties:
\begin{itemize}
\item $\bigcup_{x\in L} x+T(Q) = \ringlinsub(Q)$,
\item $(x+T(Q))\cap (y+T(Q))=\varnothing,$ for $x,y\in L$ with $x\neq y$ and
\item every intersection of $T(Q)$ with a linear subspace is measureable.
\end{itemize}
\end{Def}
Examples of fundamental domains are given in Section \ref{Sec:SpecialCases}.

The \emph{relative volume} of a subset $A\subseteq V$ denotes the volume of $A$ in the affine span $\textnormal{aff}(A)$ normalized such that any fundamental domain in $\textnormal{aff}(A)$ has volume 1. Note that it is a lower dimensional volume if the affine span of $A$ is.

We now give an explicit description of the construction of regions that determine the \ringRSs defined on full dimensional cones. In the subsequent section, we will give the whole construction in the way it occurs when starting with a polytope and determining the \ringRSs on its fcones.

\subsection{Construction on Cones}

%Fix an ambient space $V$ and a lattice $\Lambda$ in $V$ of full rank. 
We will inductively define a map $R$ from rational cones to subsets of $V$, associating a \emph{region} to each cone. From these regions, the \ringRSs can be computed via volume computations.
The author is well aware that this is a formal and very compact description of the \ringRSCc. For a step-by-step construction with examples and pictures, see Section \ref{Sec:Polytopes}. The aim of this subsection is to give a clear, short and formally precise definition of the construction. 

Let $C$ be a full dimensional rational cone in $V$. For non-full dimensional cones we intersect $V$ with the linear span of $C$ and consider that as our ambient space. For each subspace $A\subseteq V$ we  assume to have chosen a fundamental domain $T(A)$ and denote $T(B):=T(\ringlinsub(B))$ for arbitrary subsets $B\subseteq V$.

 If $C=V$ is the whole space, we set
\begin{equation*}
R(V):=T(V).
\end{equation*}

Otherwise, if $C\subsetneq V$, we assume we have constructed the regions $R(\ringfcone CF)$ for all faces $F<C$. 
Let $X_F^C$ be the set of all points $x$ in $\ringlat(\ringfcone CF)$ that fulfill the conditions:

\begin{itemize}
\item \label{Ppt:Inside} $[x+(R(\ringfcone CF)\cap \ringInt(\ringfcone CF)) ] \subseteq \ringInt(C)$\;  and
\item \label{Ppt:Nonintersect} $\left( x+ R(\ringfcone CF) \right) \cap \left(x'+ R(\ringfcone C{F'})\right) = \varnothing$ \; \hfill  \; for all $F'<C$,   with
 \raggedleft{ $F'$ \textnormal{ incomparable to } $F$  and $x'\in \ringlat(\ringfcone C{F'})$.}   
\end{itemize}

Then we define
\begin{equation*} %\label{Eqn:region}
R(C):=   \left( T(C) + \ringlinsub(C)^\perp\right)
\backslash \bigcup_{F<C} \left(X_{F}^C+R(\ringfcone CF)\right)  .
\end{equation*}

From this we can compute the values for the \emph{relative domain volume} in each region $R(C)$ for the cone $C$:
\begin{equation*} %\label{Eqn:DefDCVolume}
v_{C}:= \ringvol\big( R(C) \cap ((C \cap \Lambda ) + T ) \big).
\end{equation*}

And further the \emph{correction volumes} for each $F<C$:
\begin{equation*} %\label{Eqn:DefCorrectionVolume}
w^C_{\ringfcone CF}:=\ringvol (R(C) \cap (\ringlinsub(\ringfcone CF) \cap C))
\end{equation*}
for all faces $F<C$.

Then we get the value for $C$ as
\begin{equation*} %\label{Eqn:muIntro}
\mu(C):=v_C - \sum_{F< C} w^C_{\ringfcone CF} \cdot \mu(F).
\end{equation*}

\subsection{Construction on fcones of Polytopes} \label{Sec:Polytopes}

To give a feeling for how the construction works given an actual polytope, we go through the construction from the perspective of a general, full dimensional lattice polytope $P$. Simultaneously, we compute the values for one particular example. %To help distinguish between the general case and the example, the example will be in italicized . 
Again, for each subspace $A\subseteq V$ we choose and fix a fundamental domain $T(A)$ and set $T:=T(V)$. 

To avoid unreadable expressions, we expand $\mu, R, v$ and $w$ to functions on the faces of $P$ by setting $\mu(f):=\mu(\ringfcone P f)$, $R(f):=R(\ringfcone P f)$, $v^{f}_{g}:=v_{\ringfcone P f}$ and $w_g:=w^{\ringfcone P f}_{\ringfcone P g}$. In the example we use the same notation and keep in mind that all fcones are taken with respect to the example simplex $S$.

\hspace*{-\parindent}
\begin{minipage}{0.5\textwidth}
\paragraph{Example.} Let $S=\ringconv(v_1,v_2,v_3)\subseteq \ringR^2$ be the simplex with vertices $v_1=(1,0)$, $v_2=(2,1)$ and $v_3=(0,2)$. We consider it as a lattice polytope in $\ringR^2$ with respect to the lattice $\ringZ^2$. 

\end{minipage}
%\qquad %\qquad
\begin{minipage}{0.45\textwidth}
\qquad \qquad \scalebox{0.5}{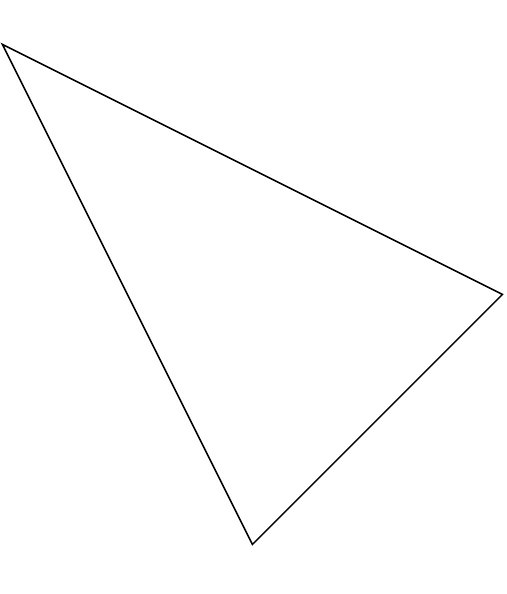} 
\end{minipage}

\vspace{12pt}
The first value we compute is $\mu(P)$ for $P$ as a face of itself. Since $\ringfcone PP=V$, we get 
\begin{equation*}
R(P)=R(V)=T
\end{equation*}
and
\begin{equation*}
v_{P}=\ringvol\big( R(P) \cap ((V \cap \Lambda ) + T ) \big) = \ringvol(T)=1, 
\end{equation*}
which directly determines the \ringRS as
\begin{equation*}
\mu(P)=v_{P}=1.
\end{equation*}
\textit{Remark.} Using this to compute the $d$-th Ehrhart coefficient, we get
\begin{equation*}
e_d =\sum_{f\in \mathcal{F}_d} \mu(f) \ringvol(f) = \mu(P)\ringvol(P) =1\cdot \ringvol(P) = \ringvol(P)
\end{equation*}
as desired, since the highest Ehrhart coefficient is known to be the relative volume of the polytope.

\paragraph{Example.}
The region $R(S)=R(\ringR^2)$ is given as the fundamental domain $T$ of $\ringZ^2$, which we choose to be the square with edgelength 1 and the origin as barycentre. Then we have $\mu(S)=v_{S}=\ringvol(T)=1.$
\\

Now, let $F<P$ be a facet of $P$. Then $H^+:=\ringfcone PF$ is a halfspace containing a hyperplane $H$ and $R(H^+)$ is defined as 
\begin{equation*}
R(H^+)= \left( T(F)+ H^\perp \right) \backslash \left( X_V^{H^+}+T  \right) ,
\end{equation*}
where 
\begin{equation*}
X_V^{H^+} = \{ x\in \Lambda \mid (x+T_V) \subseteq \ringInt(H^+) \ \}.
\end{equation*}
That means that $R(H^+)$ equals the strip $\left( T(F)+ H^\perp\right)$ minus all fundamental domains that lie entirely inside of $H^+$.

Then we have 
\begin{equation*}
v_{H^+}=\ringvol(R(H^+ )\cap ((H^+ \cap \Lambda) +T))
\end{equation*}
and 
\begin{equation*}
w_V^{H^+}=\ringvol(R(H^+) \cap  H^+).
\end{equation*}

That yields the \ringRS as 
\begin{equation*}
\mu(H^+)=v_{H^+} - w_V^{H^+} \cdot \mu(V) = v_{H^+} - w_V^{H^+}.
\end{equation*}

\begin{figure}[h]
\begin{center}
\scalebox{0.8}{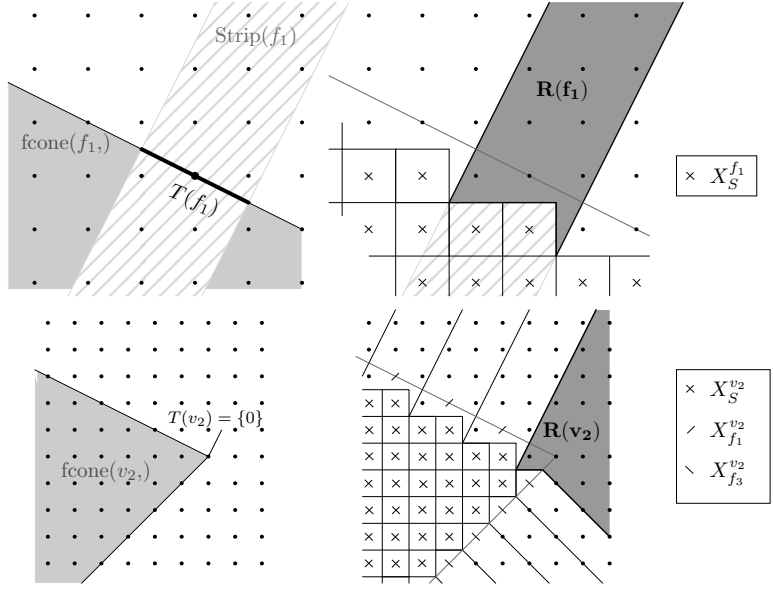} 
\end{center}
\caption{Construction of the regions $R({f_1})$ and $R({v_2})$ for the edge $f_1$ (above) and the vertex $v_2$ (below). 
}
\label{fig:ConstructRC}
\end{figure}

\begin{figure}[h] 
\begin{center}
\scalebox{0.6 }{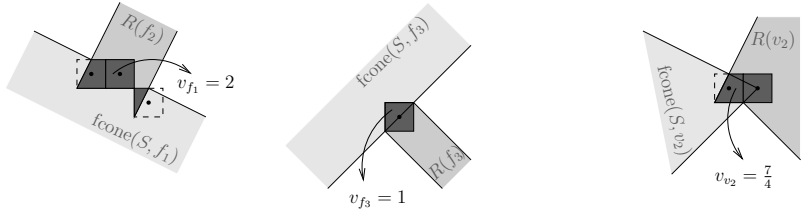}
\end{center}
\caption{Values for the relative domain volume   $v_{f}$ for certain faces $f<g$ of $S$. }
\label{fig:v}
\end{figure}

\begin{figure}[h] 
\begin{center}
\scalebox{0.6 }{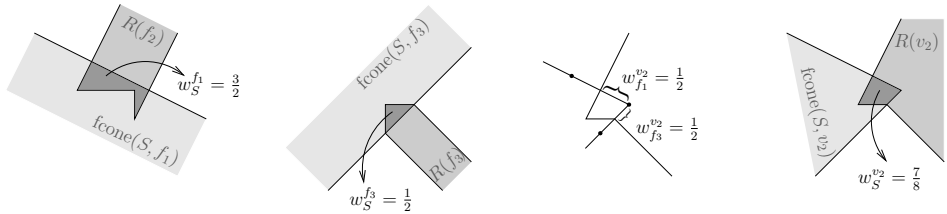}
\end{center}
\caption{Values for the \ringw $w^{f}_{g}$ for certain faces $f<g$ of $S$. }
\label{fig:w}
\end{figure}

\paragraph{Example.}
$S$ has three facets, namely the edges $f_1, f_2$ and $f_3$. Then  for $i\in \{1,2,3\}$,
$R(f_i)= \ringStrip {f_i} \backslash \left( X^{f_i}_S +T\right), $ where $\ringStrip {f_i}:=T(f_i)+\ringlinsub (\ringfcone S{f_i})^\perp$
. 
An illustration of the construction of the region $R({f_1})$ is given in Figure \ref{fig:ConstructRC}, top. The area of  $v_{{f_1}}$ is depicted in Figure \ref{fig:v} and the area of $w_{S}^{{f_1}}$ in Figure \ref{fig:w}. Altogether we get 
\vspace{-12pt}
\begin{equation*}
\mu({f_1})=  v_{{f_1}} - w_{S}^{{f_1}} = 2- 3/2 = 1/2.
\end{equation*}
Using the results from Section \ref{Sec:Codim1}, we know that the values for all facets equal $1/2$.
\\

If we now consider a codimension two face $f$ of $P$, then we have exactly two facets $F_1$ and $F_2$ of $P$ that meet in $f$. That means that $W:=\ringfcone Pf$ is a wedge defined by the intersection of the halfspaces $H_1^+:=\ringfcone P{F_1}$ and $H_2^+:=\ringfcone P{F_2}$, whose lineality spaces are  the hyperplanes $H_1:=\ringlinsub(\ringfcone P{F_1})$ and $H_2:=\ringlinsub(\ringfcone P{F_2})$, respectively. Hence, the lineality space of $\ringfcone Pf$ is the line $ \ringlinsub(\ringfcone Pf)=H_1\cap H_2$.

As above, for $\ringfcone PP=V$ we have 
\begin{equation*}
X_V^W=\{ x\in \Lambda \mid (x+T_V) \subseteq \ringInt(W) \ \},
\end{equation*}
 for $H_1^+$ we have
\begin{align*}
X_{H_1^+}^W = \{ x \in \ringlat(H_1^+) \mid & (x+R(H_1^+)\cap \ringInt(H_1^+)) \subseteq \ringInt(W) \textnormal{ and } \\ 
&(x+R(H_1^+))\cap (x'+R(H_2^+))=\varnothing \textnormal{ for all } x'\in \ringlat(H_2^+)\}
\end{align*}
and analogously for $H_2^+$
\begin{align*}
X_{H_2^+}^W = \{ x \in \ringlat(H_2^+) \mid & (x+R(H_2^+)\cap \ringInt(H_2^+)) \subseteq \ringInt(W) \textnormal{ and } \\ 
&(x+R(H_2^+))\cap (x'+R(H_1^+))=\varnothing \textnormal{ for all } x'\in \ringlat(H_1^+)\}.
\end{align*}

Then the region $R(W)$ is given by
\begin{align*}
R(W)= &\left( T(W) + \ringlinsub(W)^\perp\right)\\ 
		&\backslash \left( (X_{V}^W+T ) \cup (X_{H_1^+}^W+R(H_1^+) )\cup (X_{H_2^+}^W+R(H_2^+) ) \right)  
	.
\end{align*}

The relative domain volume is 
\begin{equation*}
v_{W}=\ringvol(R(W)\cap ((W \cap \Lambda) +T))
\end{equation*}
and we have to consider three correction volumes:
\begin{align*}
w^W_{H_1^+}&=\ringvol(R(W) \cap (H_1\cap C)), \\
w^W_{H_2^+}&=\ringvol(R(W) \cap (H_2\cap C)) \textnormal{ and }\\
w_V^{W}&=\ringvol(R(W) \cap W).
\end{align*}
The \ringRS for $W$ then is
\begin{align*}
\mu(W)= v_{W} - w^W_{H_1^+} \cdot \mu(H_1^+)  - w^W_{H_2^+} \cdot \mu(H_2^+)  - w^W_{V} \cdot \mu(V) 
\end{align*}

\paragraph{Example.}
To compute the value corresponding to the vertex $v_2$ of $S$, we construct the region $$R({v_2})=\ringR^2 \backslash \left( (X^{v_2}_S + T) \cup (X^{v_2}_{f_1} + T(f_1)) \cup(X^{v_2}_{f_3} + T(f_3))  \right)$$ 
as shown in Figure \ref{fig:ConstructRC}, bottom. We then get $v_{{v_2}}=7/4$ (cf. Figure \ref{fig:v}) and as correction volumes we have the 1-dim. relative volumes $w^{{v_2}}_{{f_1}}=1/2$ and $w^{{v_2}}_{{f_2}}=1/2$ and the 2-dim. relative volume $w^{{v_2}}_{S}=7/8$ (cf. Figure \ref{fig:w}). 
Altogether we get the \ringRS for the fcone of $S$ at the vertex $v_2$ as
\begin{align*}
\mu({v_2}) & = v_{{v_2}} - w^{{v_2}}_{{f_1}} \cdot \mu({f_1})
	- w^{{v_2}}_{{f_3}} \cdot \mu({f_3})
	- w^{{v_2}}_{{S}} \cdot \mu({S}) \\
	& = \frac{7}{4} - \frac{1}{2} \cdot \frac{1}{2} - \frac{1}{2} \cdot \frac{1}{2} - \frac{7}{8} \cdot 1 = \frac{3}{8}.
\end{align*}
Analogously, we can compute $\mu({v_3})=1/4$ and $\mu({v_1})=3/8$. For reasons of symmetry, the latter has to equal the value $\mu({v_2})$, see Section \ref{Sec:Symmetry} for an explanation.

For the general faces $f$ of $P$ the region of $f$ is computed as 

\begin{equation*} 
R(f):=\left( T(f) + \ringlinsub(\ringfcone P f)^\perp\right) 
\backslash \bigcup_{g>f} \left(X_{g}^{f}+R(g)\right) ,
\end{equation*}

the relative domain volume is
\begin{equation*} 
v_{f}:= \ringvol\big( R(f) \cap ((f \cap \Lambda ) + T ) \big)
\end{equation*}

and the correction volume for each face $g$ of $P$ with $g>f$ is
\begin{equation*} 
w^{f}_{g}:=\ringvol (R(f) \cap (\ringlinsub(\ringfcone Pg) \cap f)).
\end{equation*}

Then we get the value for $\ringfcone Pf$ as
\begin{equation*} 
\mu(f):=v_{f} - \sum_{g>f} w^C_{g} \cdot \mu(g).
\end{equation*}

\section{The geometrical connection between the \ringRSs and the number of lattice points}\label{Sec:Motivation}

This section gives a motivation behind the definition of the \ringRSs and, along the way, an outline of the proof. Complete proofs can be found in \cite{RS17}.

Let $P=\ringconv(v_1,\ldots,v_m)$ be a full dimensional lattice polytope in $d$-dimensional Euclidean space $V$. We use the same notation as in Section \ref{Sec:Polytopes}, where we write $\mu(f)$ for $\mu(\ringfcone Pf)$ and the same for $R$, $v$ and $w$.

For a face $f$ of $P$ and $t\in \ringZgeq$ define the set $\ringX(tf)$ of all \emph{feasible lattice points} in $t\cdot f$ as
\begin{equation*}
\ringX(tf):= \bigcap_{i=1}^m \left( X^{{v_i}}_{{f}} +tv_i \right).
\end{equation*}

After the construction of all regions for $P$, we get the following tiling of $V$
\begin{thm}[{\cite[Thmeorem 1]{RS17}}]
 \label{Thm:Tiling}
Let $P\subseteq V$ be a full dimensional 
lattice polytope. There exists a $t_0\in \mathbb{Z}_{>0}$ such that for each $t\geq t_0$   we have a tiling of $V$
into translated regions of the form
\begin{equation*}
\{x+R(f) : \; f\leq P, \, x\in \ringX(tf)\}.
\end{equation*}
\end{thm}

\begin{figure}[h] 
\begin{center}
\scalebox{0.269 }{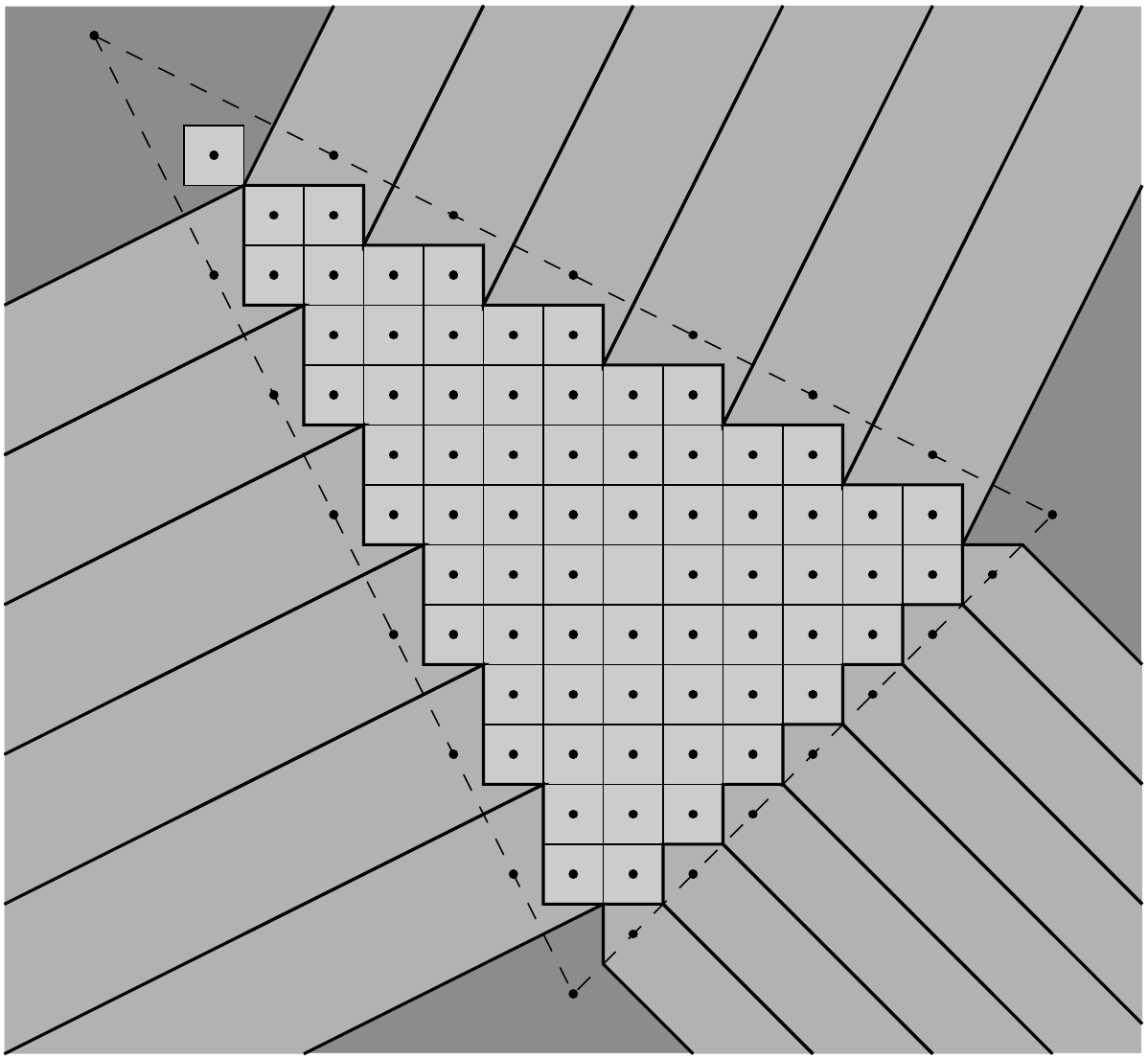} \!   % \qquad% \qquad 
\scalebox{0.269 }{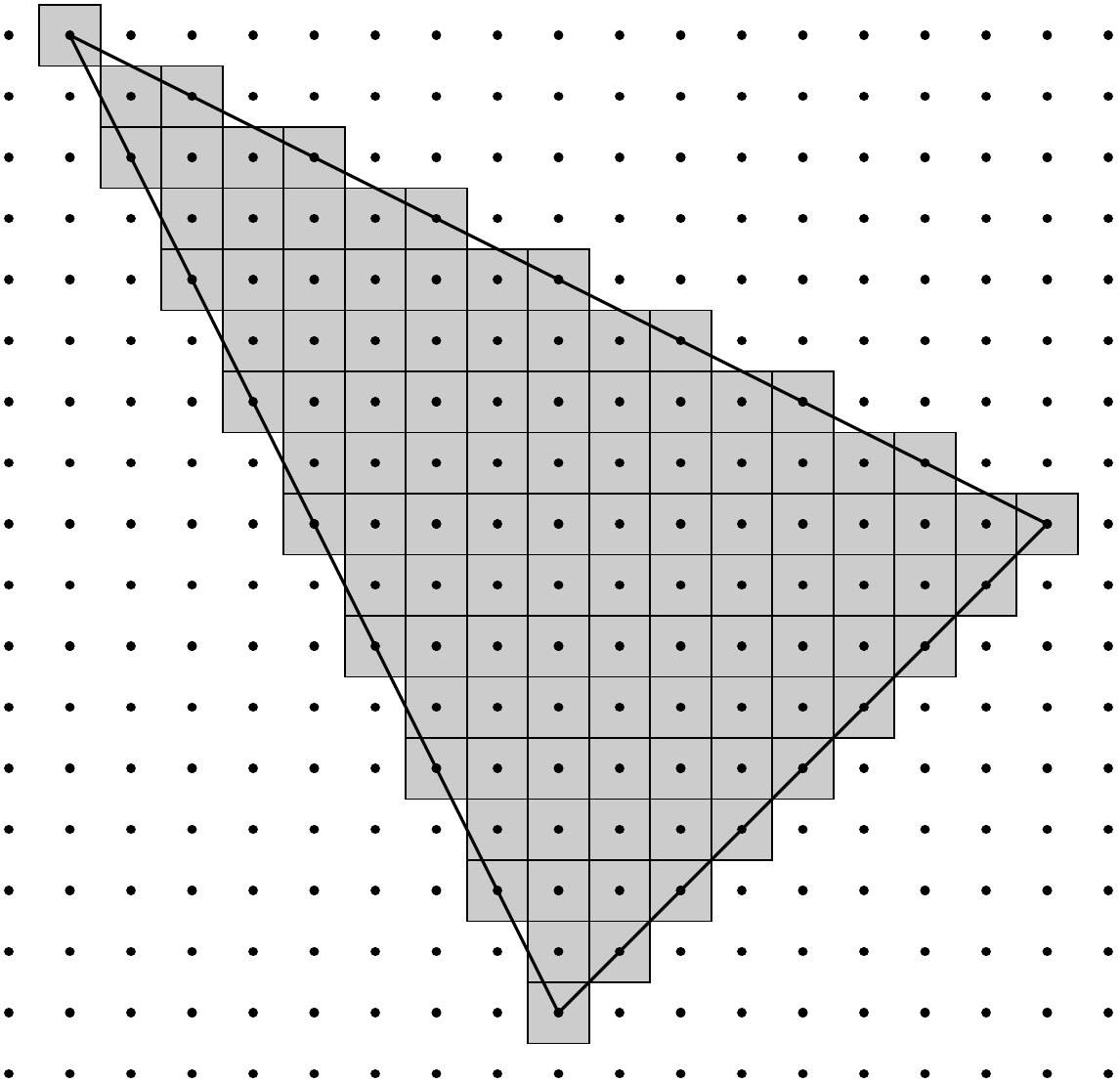} \! 
\scalebox{0.269}{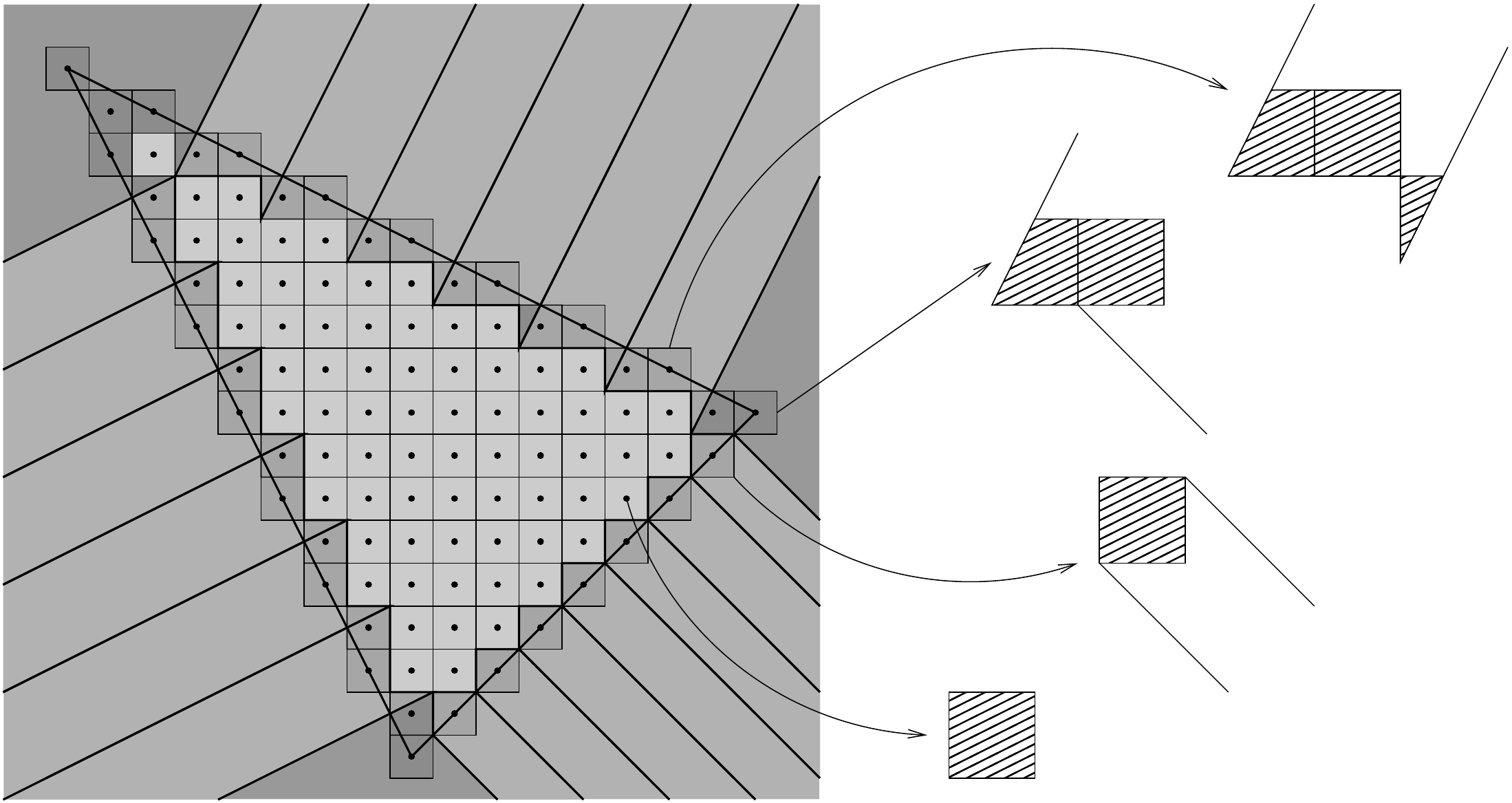}
\end{center}
\caption{ \textbf{Left:} Tiling of $\ringR^2$ from Thm. \ref{Thm:Tiling} for the triangle $S$ with $t=8$; \textbf{middle:} The (fundamental) domain complex of $8 S$; \textbf{right:} Values of $v_f$ for some faces $f<S$}
\label{fig:Tiling&DC}
\end{figure}

\paragraph*{Example} The tiling of $8\cdot S$ for the example polygon $S$ from Section \ref{Sec:Polytopes} can be seen in Figure \ref{fig:Tiling&DC} on the left; in the picture, each region is translated by a unique lattice point marked $\bullet$ and the union of  all these  lattice points within a certain face $f$ form the set $\ringX(8f)$.

%The set $\ringX(tf)$ of all feasible lattice points in $tf$ might be empty for small $t$, such that the condition $t\geq t_0$ is  necessary. 

The crucial step between counting the number of lattice points in $P$ and the \ringRS given by volumes is to interpret the number of lattice points in $tP$ as the volume of all translates of  $T$ by the lattice points in $P$:
\begin{equation*} 
|\Lambda \cap tP| = \sum_{x\in \Lambda\cap tP} \ringvol(x+T) = \ringvol\underbrace{((\Lambda\cap tP)+T)}_{=:DC},
\end{equation*}
The first part of the equation holds, since by definition $\ringvol(T)=1$ for any fundamental domain~$T$ of~$\Lambda$, and the second part follows from $(x+T)\cap (y+T)=\varnothing$ for all $x,y\in \Lambda$ with $x\neq y$. 
We call the set $DC:=(\Lambda \cap tP)+T$ 
a \emph{(fundamental) domain complex} of $tP$ (cf. Figure
\ref{fig:Tiling&DC}, middle). By taking the volume of the respective
part of the domain complex in each region of the tiling in
Theorem~\ref{Thm:Tiling}, we get
\begin{equation*}
|\Lambda \cap tP|= \ringvol(DC) = \sum_{f\leq P} \sum_{x\in \ringX(tf)}\underbrace{ \ringvol \big( (x+R(f)) \cap DC \big)}_{(*)}.
\end{equation*}

It turns out that $(*)$ is exactly the value $v_{f}$ as defined in Section \ref{Sec:Construction}, which thus yields
\begin{equation} \label{Eqn:SumLatticepoints}
|\Lambda\cap tP|=  \sum_{f\leq P} \sum_{x\in \ringX(tf)} v_{f} =  \sum_{f\leq P}  v_{f} \cdot |\ringX(tf)|,
\end{equation}
see Figure \ref{fig:Tiling&DC}, right.

\paragraph*{Example} Equation \eqref{Eqn:SumLatticepoints} for $4S$ gives 
\begin{align*}
|\ringZ^2\cap 4S|=109 &= 1\cdot 70 + 2\cdot 6 + 2 \cdot 6 + 1 \cdot 7 + \frac{7}{4} \cdot 1 + \frac{7}{4} \cdot 1 + \frac{9}{2} \cdot 1 \\
	&= \sum_{f\leq S}  v_{f}\cdot |\ringX(8f)|,
\end{align*}
where the sum runs over the faces $S, f_1, f_2, f_3, v_1, v_2, v_3$ in that order.
\\

To extract the value for the local formula from Equation \eqref{Eqn:SumLatticepoints}, we need to determine the difference between $|\ringX(tf)|$ and $\ringvol(f)$. 
It can be shown (see \cite[Lemma 8]{RS17}) that by the construction of the correction volumes $w$ we have that 
\begin{equation} \label{Eqn:vol(tf)}
\ringvol(tf)= \sum_{g\leq f}  w^g_f \cdot |\ringX(tf)| ,
\end{equation}
where we formally define $w^f_f=1$. 

We  want the \ringRS to be defined such that 
\begin{align*}
|\Lambda \cap t P|
		& = \sum_{f\leq P} \mu(f)\cdot \ringvol(f) 
\end{align*}
holds.
Using Equation (\ref{Eqn:vol(tf)}), we get
\begin{align*}
\sum_{f\leq P} \mu(f)\cdot \ringvol(f) 
		& =
\sum_{f\leq P} \mu(f)\cdot \left( \sum_{g\leq f}  w^g_f \cdot |\ringX(tf)| \right).
\end{align*}
By combinatorially reordering the right hand side, we get
\begin{align*}
\sum_{f\leq P} \mu(f)\cdot \ringvol(f) 
		& = \sum_{f\leq P} \left[ \mu(f)\cdot |\ringX(tf)| + \mu(f) \cdot \sum_{ h> f} w^f_h \cdot |\ringX(tf)|   \right]\\
		& = \sum_{f\leq P} \left[  \left( \mu(f) + \sum_{h>f} \cdot w^f_h \right)\cdot |\ringX(tf)|  \right]
\end{align*}

Hence, we want to define the \ringRS such that the following equation holds
\begin{equation*}
\sum_{f\leq P}  v_{f} \cdot |\ringX(tf)| = |\Lambda\cap tP| 
		=\sum_{f\leq P} \left[  \left( \mu(f) + \sum_{h>f} \cdot w^f_h \right)\cdot |\ringX(tf)|  \right].
\end{equation*}
Comparing coefficients we get
\begin{align*}
		v_{f} & =\mu(f) + \sum_{h>f} \cdot w^f_h,
\end{align*}
which directly leads to the definition of the \ringRSs as
\begin{align*}
		\mu(f) & = 	v_{f} - \sum_{h>f} \cdot w^f_h.
\end{align*}

%((((Sideremark: Let $\ringStrip F$ be the strip $\left( T(F)+ H^\perp\right)$.  Another way to compute $v_{H^+}$ is to compute the volume of all  the fundamental domains that intersect $H$ and whose lattice point lies in the strip and in $H^+$. That is:
%\begin{align*}
%v_{H^+}&=\ringvol\left( \left\{  x \in \Lambda\cap \ringStrip F \cap H^+ \mid (x+T)\cap H \neq \varnothing \right\} +T \right) \\
%	&=|\left\{  x \in \Lambda\cap \ringStrip F \cap H^+ \mid (x+T)\cap H \neq \varnothing \right\}|
%\end{align*}))))

\section{Symmetry and special cases} \label{Sec:SpecialCases}

Symmetry of the \ringRSs can be achieved by choosing symmetric fundamental domains, for example by taking  \emph{Dirichlet--Voronoi cells}, as we will show in Section \ref{Sec:Symmetry}. We will finish this article by proving that given central symmetry, the values on halfspaces are determined to be $1/2$. That implies that the value  is always $1/2$ for fcones of a facet of a polytope.

\subsection{Dirichlet--Voronoi cells and symmetry of RS-$\mu$} \label{Sec:Symmetry}

Possibly the most natural choice of fundamental domains are \emph{Dirichlet--Voronoi cells}. %For lattices sometimes referred to as 'Wigner-Seitz cells' 
Given a space $V$ and an inner product $\langle \cdot  ,\cdot \rangle$ with induced norm $\| \cdot \|$, the Dirichlet--Voronoi cell of a sublattice $L \subseteq \Lambda$ is defined as 
\begin{equation*}
\textnormal{DV}(L, \langle\cdot,\cdot\rangle):=\{x\in \ringlin(L) : \|x\| \leq \|x-a\| \textnormal{ for all } a\in L\}. 
\end{equation*}
In this definition, it is not yet a fundamental domain of the lattice~$L$, but by considering  the Dirichlet--Voronoi cell half open, it can be seen as a fundamental domain of the lattice.
Dirichlet--Voronoi cells are naturally centrally symmetric and can be forced to have certain symmetries by choosing a suitable inner product.

Let $P$ be a lattice polytope and $\ringG$ a subgroup of all lattice symmetries of $P$, i.e. $\ringG$ is a finite matrix group with $A\cdot P:=\{A\cdot x: x\in P\}=P$ and $A\cdot \Lambda = \Lambda$ for all $A\in \ringG$. Then we can define a $\ringG$-invariant inner product by taking 
\begin{align}
\langle x,y\rangle_\ringG :=x^t G y & & \textnormal{ for all } x,y\in V,
\end{align}
with the Gram matrix  $G$ given by
\begin{equation}
G:= \frac{1}{|\ringG|} \sum_{A\in \ringG} A^tA.
\end{equation}
Let $\|\cdot\|_\ringG$ be the induced norm and
let $D$ be the Dirichlet--Voronoi cell for $\Lambda$ given by the that particular inner product,
\begin{equation*}
D:=\textnormal{DV}(\Lambda, \langle\cdot,\cdot\rangle_\ringG)=\{x\in V : \|x\|_\ringG \leq \|x-p\|_\ringG \textnormal{ for all } p\in \Lambda\}. 
\end{equation*}

Then $D$ is invariant under the action of $\ringG$:
Let $x\in D$, then for $A\in \ringG$ we have
\begin{align*}
\|Ax\|_\ringG = \|x\|_\ringG \leq \| x- p\|_\ringG = \| A x - A p\|_\ringG & &\textnormal{ for all } p \in \Lambda.
\end{align*}
Since $ A \Lambda = \Lambda$, we get $AD\subseteq D$ for all $A\in \ringG$. Substituting $A$ by $A^{-1}$, we get $A^{-1} D\subseteq  D $ which yields $ D\subseteq AD$ and hence $AD=D$.
Similarly, we see that for all faces $f$ in the same $\ringG$-orbit the
fcones and Dirichlet--Voronoi cells in $\Lambda \cap N_f^\perp$
are mapped onto each other. Hence, the used regions are invariant under the action of $\ringG$ and $\mu$ is constant on $\ringG$-orbits.

\subsection{Codimension one faces under central symmetry} \label{Sec:Codim1}

It is known that the second highest Ehrhart coefficient always 
equals $1/2$ times the sum over the relative volumes of the 
facets of a polytope. A natural conjecture would be that all 
values of McMullen's formulas corresponding to facets 
(in this case all values on halfspaces) have the value $1/2$. 
This is not true in general for the \ringRSsc, but we show here 
that it does hold for the \ringRSs when all fundamental domains 
are centrally symmetric.
Therefore, in the following let $T(A)$ be a centrally symmetric fundamental domain for each $A\subseteq V$. Again, we denote $T:=T(V)$. 

\begin{figure}[h] 
\begin{center}
\scalebox{0.3}{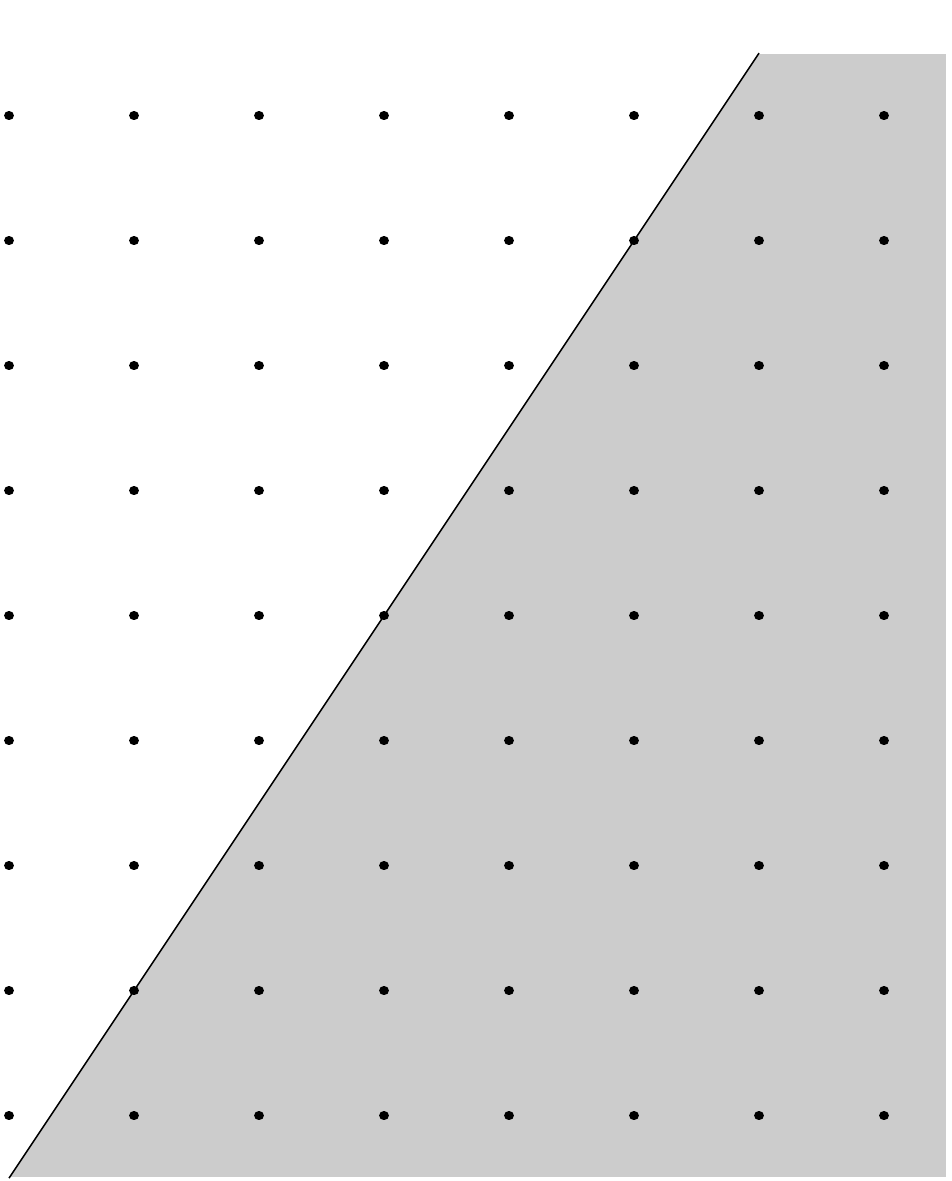} \qquad
\scalebox{0.3}{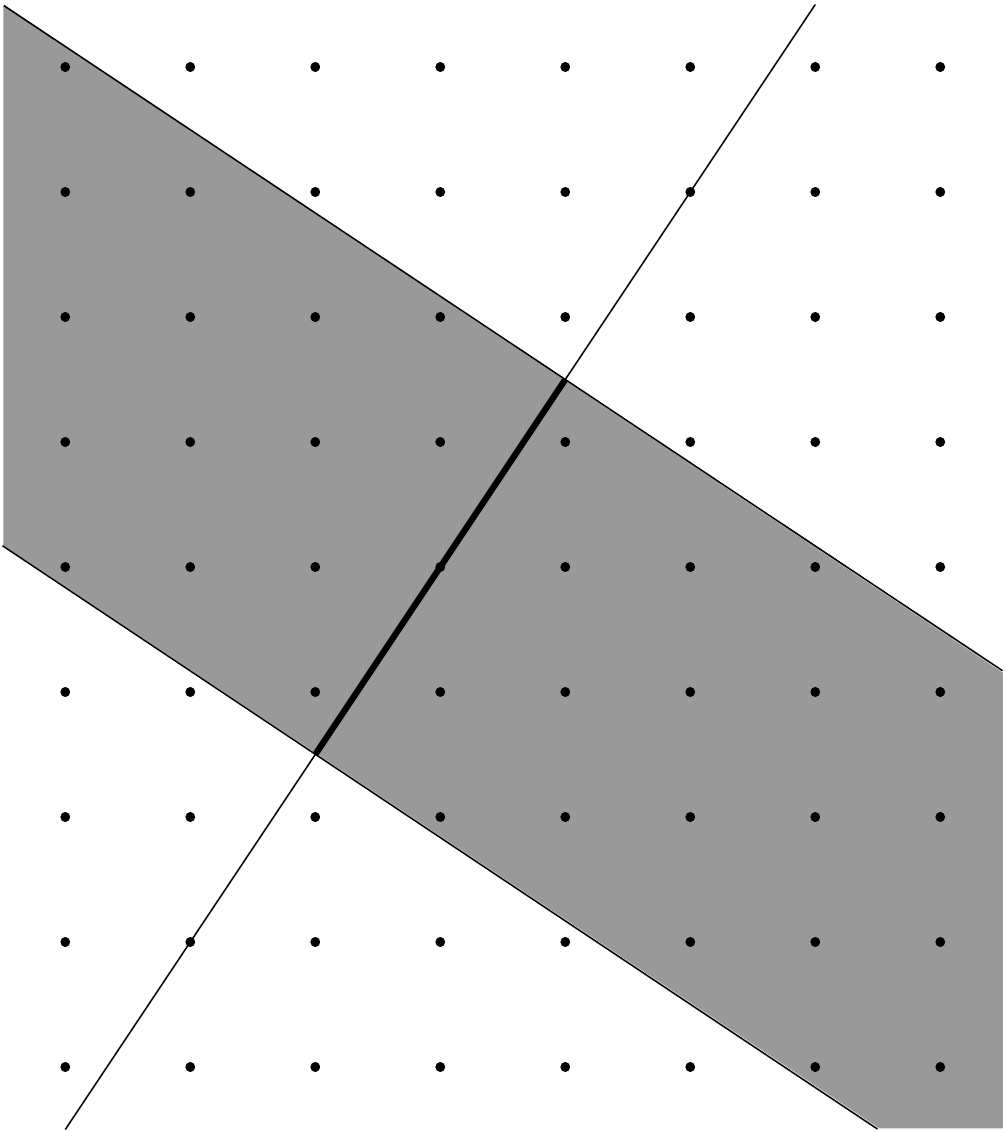}
\end{center}
\caption{\textbf{left}: $C$ with $C^+$, $C^-$ and $\ringlinsub(C)$; \textbf{right}: $\ringStrip C=T(C)+\ringlin(C)^\perp$}
\label{fig:strip}
\end{figure}

Now, let $F$ be a facet of $P$ with fcone $C:=\ringfcone PF$. That means $\ringlinsub(C)$ is a hyperplane in $V$. Let $C^+$ be the open halfspace inside of $C$ and $C^-$ the open halfspace on the other side, i.e. the complement of $C$. 

Then 
\begin{align*}
R(C)= \underbrace{( T(C) + \ringlin(C)^\perp)}_{=:\ringStrip C} \backslash \left( X_V^{C}+T  \right) \\
= \ringStrip C \cap \left( (\Lambda \backslash  X_V^{C})+T  \right),
\end{align*}
where 
\begin{equation*}
X_V^{C} = \{ x\in \Lambda \mid (x+T) \subseteq C^+ \ \}.
\end{equation*}

\begin{figure}[h] 
\begin{center}
\scalebox{0.3}{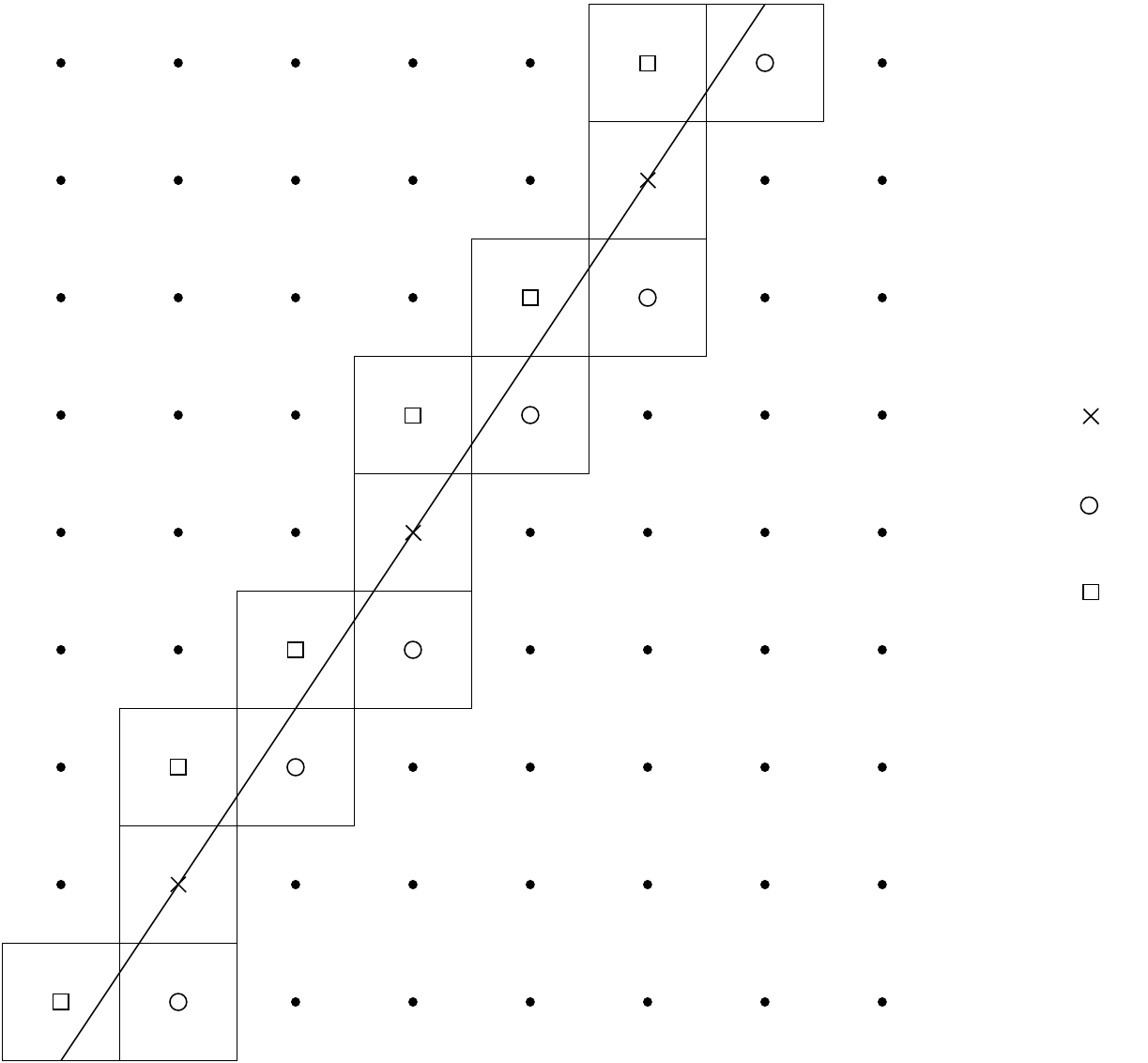} \quad
\scalebox{0.3}{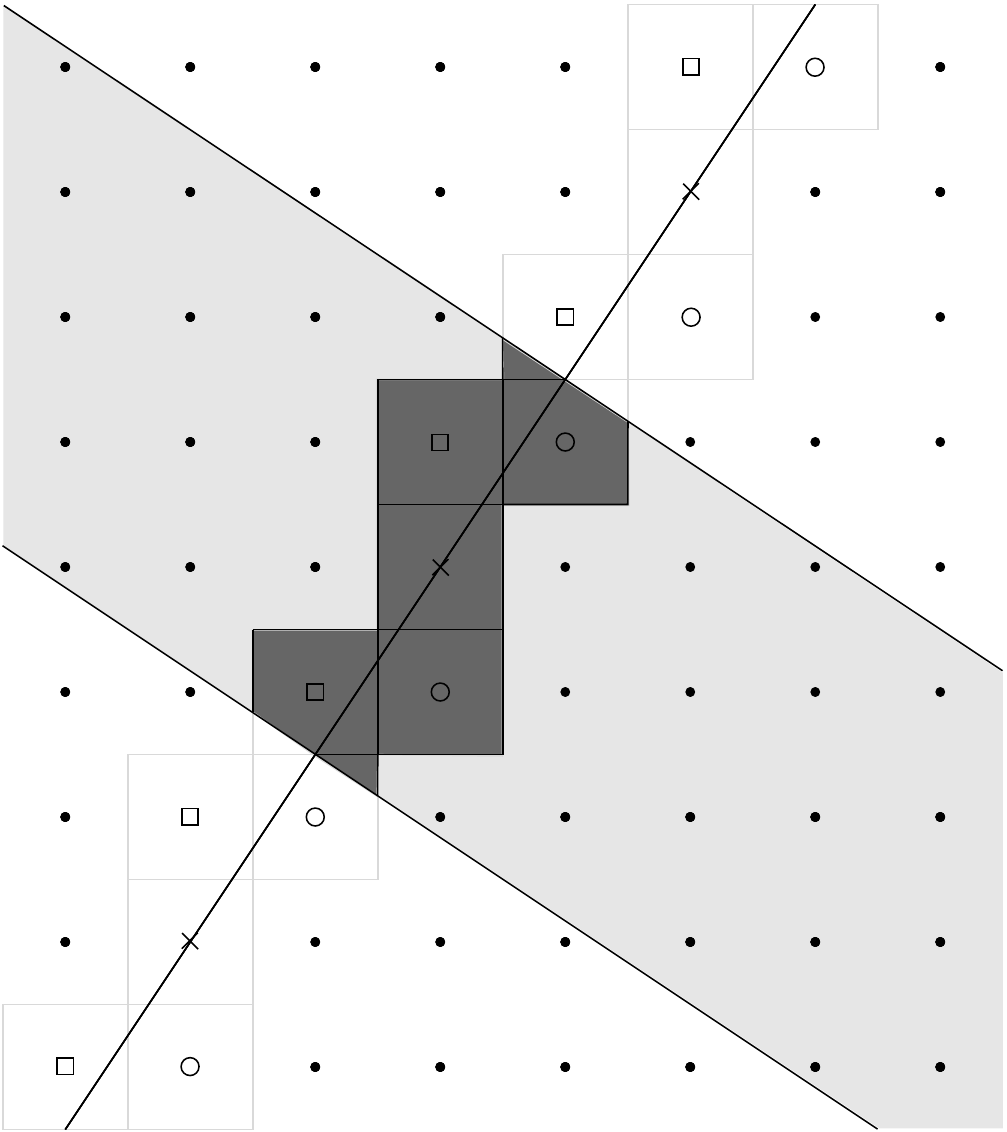}
\end{center}
\caption{$(X+T)$ and $\ringReg(C)$}
\label{fig:X}
\end{figure}

Taking a closer look at the construction in Section \ref{Sec:Construction}, we note that for computing $v_C$ and $w^C_K$ and hence $\mu(C)$ it is only necessary to consider $R(C)$ intersected with the union of all fundamental domains that have a nonempty intersection with $\ringlinsub(C)$. 
We therefore consider the \emph{relevant part} $\ringReg(C)$ of the region $R(C)$:
\begin{align*}
\ringReg(C) 	=&(\underbrace{\{  p\in \ringlat(W) \mid  (p+T)\cap \ringlinsub(C) \neq \varnothing \}}_{:=X}+T)  \cap \ringStrip C
\end{align*}
This relevent part of $R(C)$ can be partitioned into three parts
\begin{align*}
\ringReg(C)	 	=& (X + T) \cap \ringStrip C \\
		=& (\underbrace{(X\cap \ringlinsub(C)}_{:=X_0}+T)\cap \ringStrip C\\ 
		   & \cup   (\underbrace{(X\cap C^+}_{:=X_+}+T)\cap \ringStrip C\\ 
		   & \cup  (\underbrace{(X\cap C^-}_{:=X_-}+T)\cap \ringStrip C,
\end{align*}
where the unions are disjoint, since $X_0, X_+, X_-$ are. For an illustration see Figure~\ref{fig:X}.

Recall that $\mu(C)$ is defined as 
\begin{align*}
\mu(C)	&= v_C - \mu(W) \cdot w^C_W \\
		&= \ringvol(\ringReg(C) \cap ((C\cap \ringLat(W)))+T) - 1\cdot \ringvol(\ringReg(C) \cap (W \cap C))
\end{align*}

To show that  $\mu(C)=1/2$, we show that everything but half the lattice cell around the origin cancels out nicely. We use the fact that everything is centrally symmetric in the following sense:

Let $\sigma_0$ be the point reflection at the origin:
\begin{align*}
\ringsymm \colon 	& V \rightarrow V \\
			& v \mapsto -v
\end{align*}

Then $\ringsymm(T)=T$ by assumption and since $C$ is a halfspace, we also have
\begin{align*}
\ringsymm (\ringlinsub(C))		&=\ringlinsub(C) ,\\
\ringsymm(C^+)		&=C^-, \\ 
\ringsymm(C^-)		&=C^+ ,\\
%\ringsymm(\FundCell(C))	&=\FundCell(C),\\ 
\ringsymm(\ringStrip C)		&=\ringStrip C\\
\ringsymm(\ringReg(C))	&=\ringReg(C)
\end{align*}

Now observe that $\ringvol((X_0+T)\cap \ringStrip C) =\ringvol( 0+T)=1$, since the intersection of a centrally symmetric fundamental domain with a linear subspace is always contained in the centrally symmetric fundamental domain of that subspace. 

 %\note{fundamental domain for strip can only get bigger but not as big such that $p+T$ can be included for any other $p\in \Lambda\backslash \{0\}$ ...}

Since $\ringsymm$ does not change the volume, we can use it to gain information on the occuring volumes:
\begin{align*}
\ringsymm((T\cap \ringStrip C ) \cap C^+)=T\cap \ringStrip C  \cap C^-
\end{align*}
and hence, 
\vspace{-12pt}
\begin{align*} \vspace{-12pt}
\ringvol(((X_0+T)\cap \ringStrip C)\cap C^-) =\frac{1}{2}.
\end{align*}

\vspace{-12pt}
\begin{figure}[h] 
\begin{center}
\scalebox{0.25}{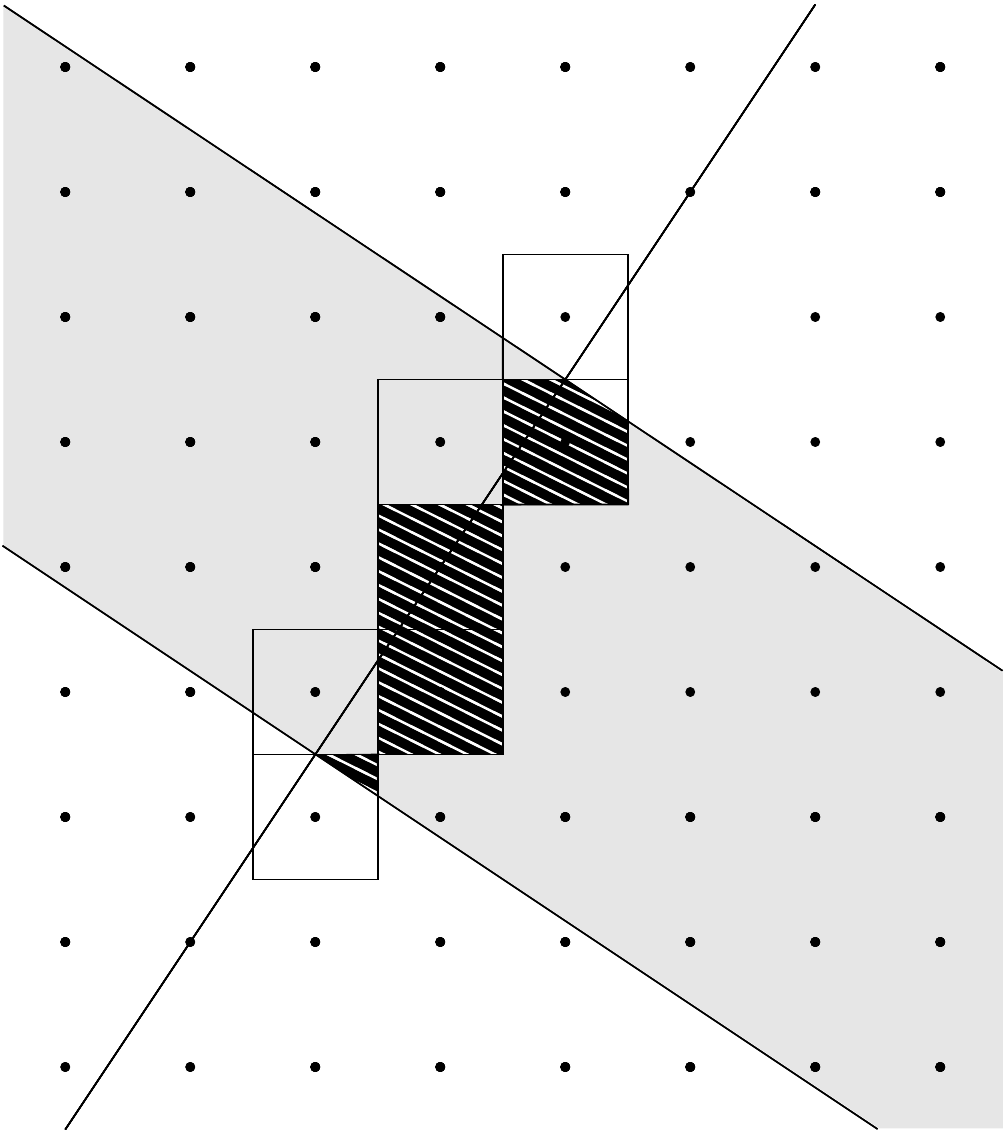} \qquad
\scalebox{0.25}{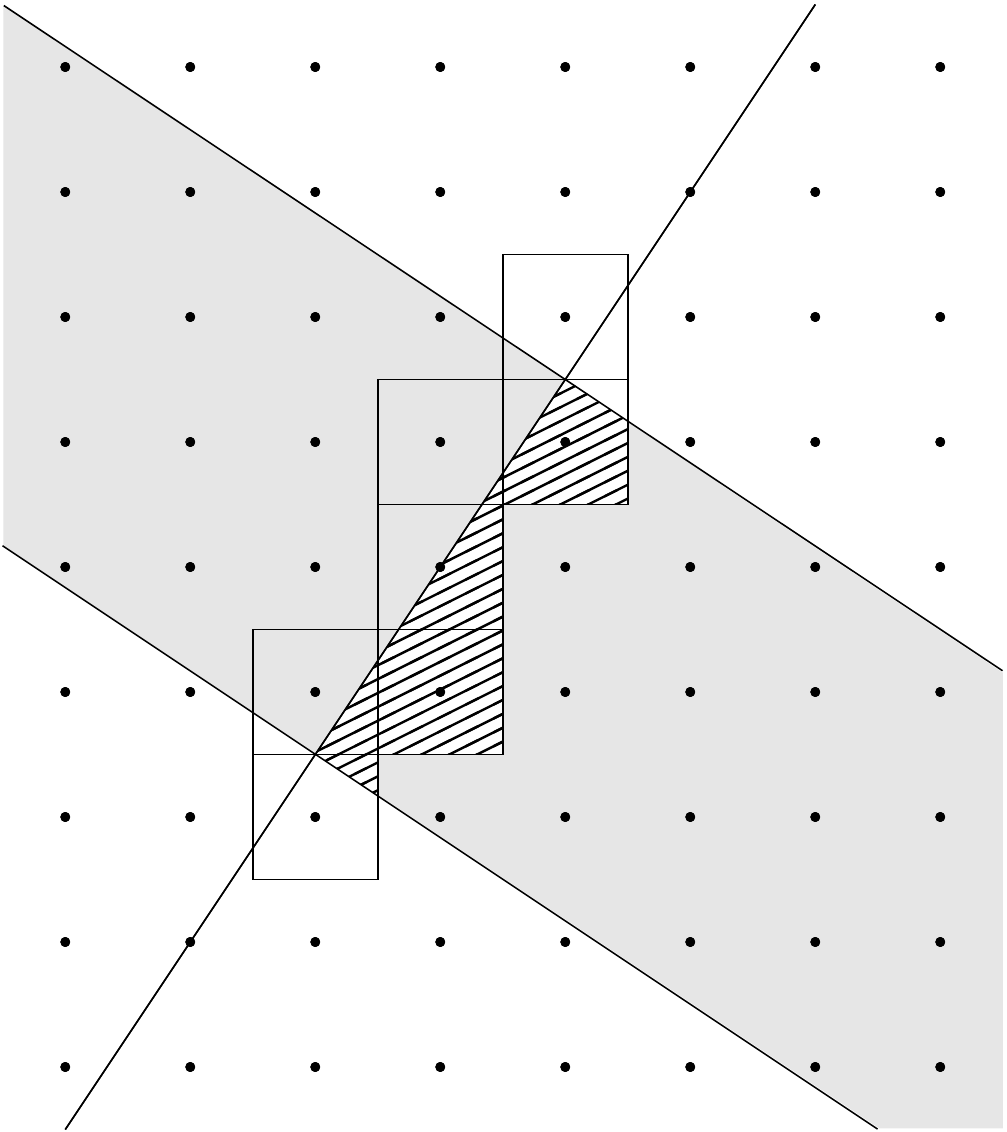} \hspace{-5pt}
\scalebox{0.25}{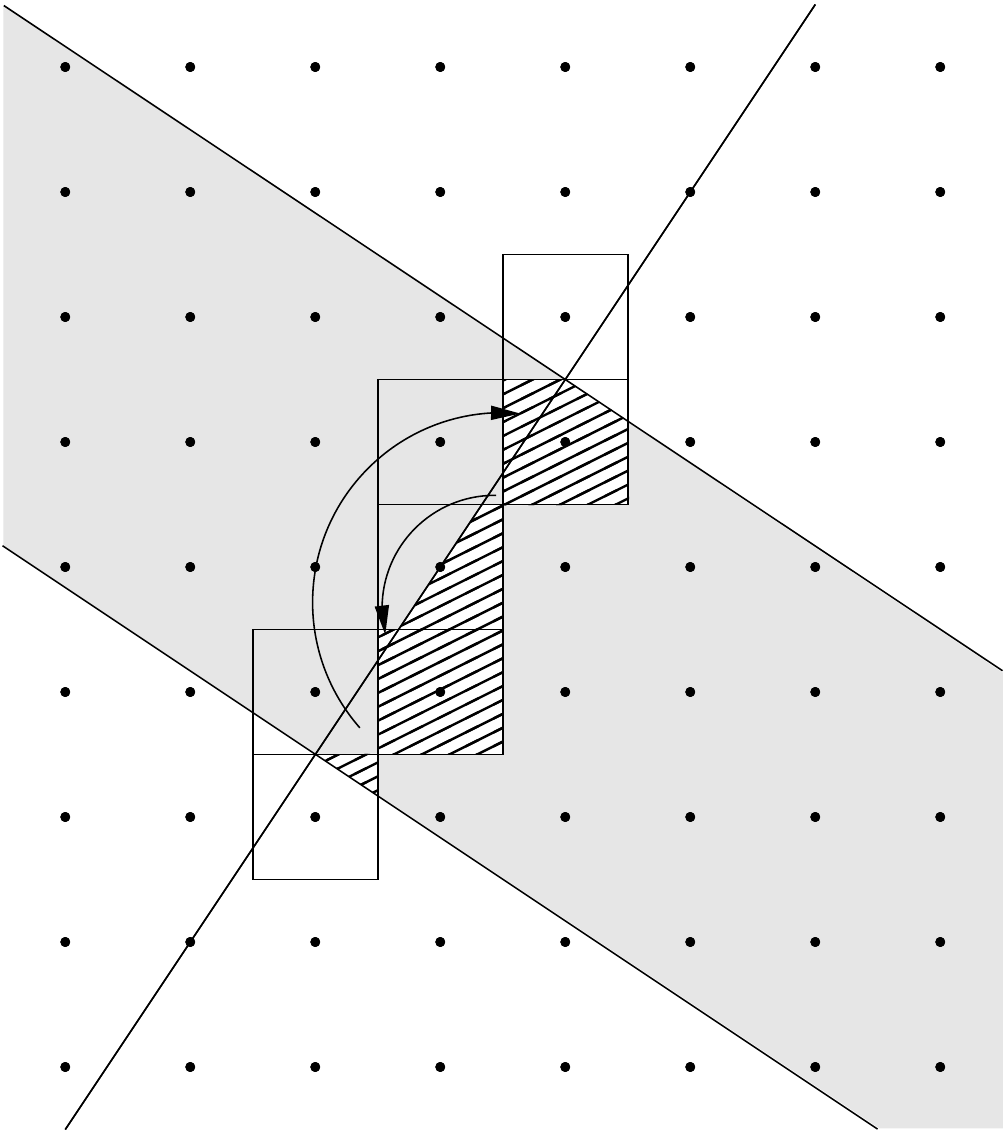} \qquad
\scalebox{0.25}{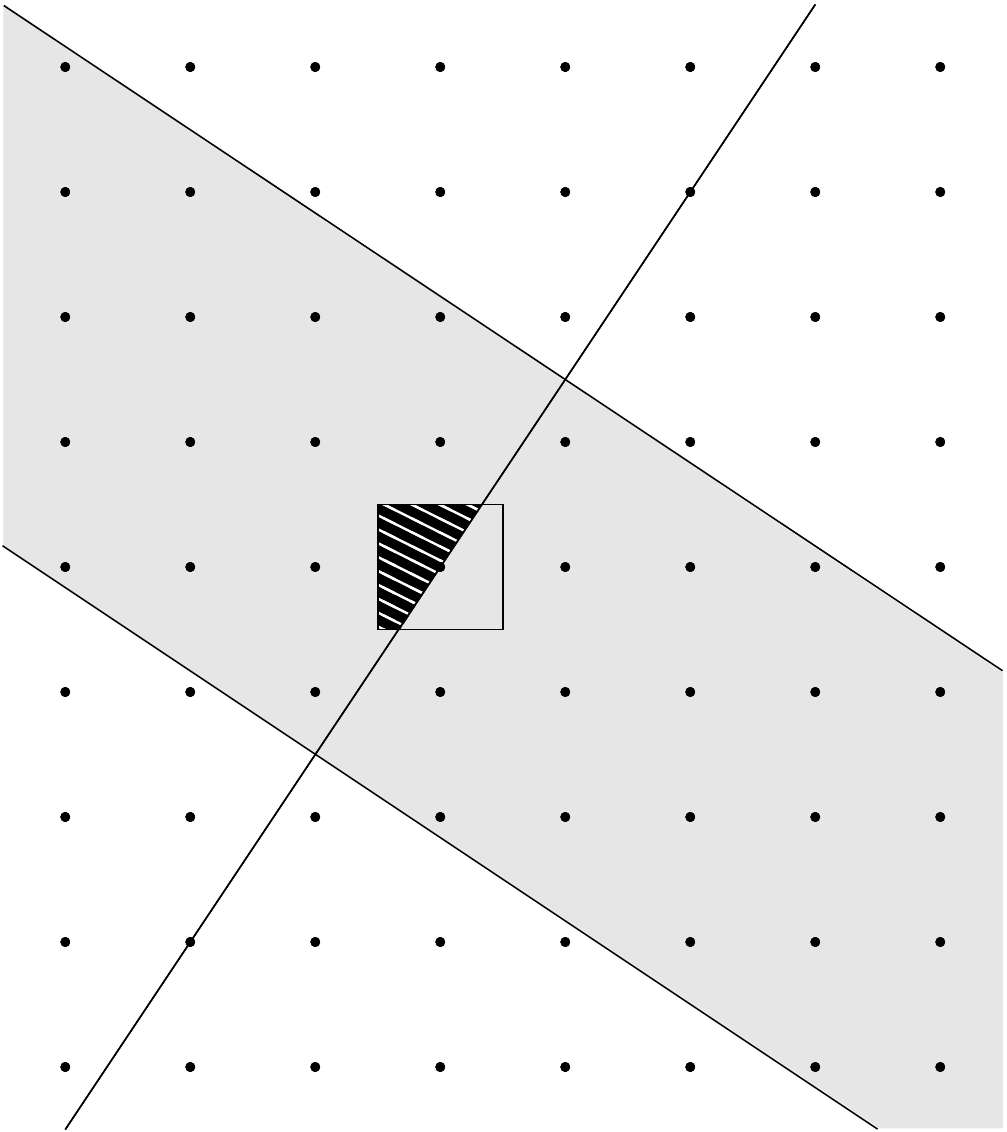}
\end{center}
\caption{ \textbf{Left:} $v_C$, \textbf{middle:} $w^C_W$ and  \textbf{right:} $w^C_W-v_C$}
\label{fig:vw}
\end{figure}

Moreover, we have the disjoint union:
\begin{align*}
\ringReg(C) \cap   C^+ = 	& ((X_0+T)\cap \ringStrip C) \cap C^+ \\ 
					& \cup (X_+ +T)\cap \ringStrip C) \cap  C^+ \\
					& \cup  (X_- +T)\cap \ringStrip C) \cap  C^+ 
\end{align*}
and since  $\ringsymm(((X_- +T)\cap \ringStrip C) \cap  C^+ )= ((X_+ +T)\cap \ringStrip C) \cap  C^-$, the  two have equal volume and we have 
 
\begin{align*}
w^C_V  = &\ringvol(\ringReg(C) \cap (V \cap C^+)) \\
  = & \ringvol( ((X_0+T)\cap \ringStrip C) \cap C^+ ) \\ 
					 &+\ringvol( (X_+ +T)\cap \ringStrip C) \cap  C^+) \\
					 &+\ringvol(  (X_+ +T)\cap \ringStrip C) \cap  C^- ) \\
				=&	\ringvol( ((X_0+T)\cap \ringStrip C) \cap C^+ ) \\ 
					&+\ringvol( (X_+ +T)\cap \ringStrip C) ),
\end{align*}
see Figure \ref{fig:vw}, middle.                      
Together with
\begin{equation*}
v_C= \ringvol(\ringReg(C) \cap ((C\cap \Lambda)+T)) = \ringvol((X_0+T)\cap \ringStrip C) +\ringvol((X_++T)\cap \ringStrip C), 
\end{equation*}
(cf. Figure \ref{fig:vw}, left), we finally get
\begin{align*}
\mu(C)&= v_C - w^C_V \\
& = \ringvol(\ringReg(C) \cap ((C\cap \Lambda))+T) -  \ringvol(\ringReg(C) \cap (V \cap C)) \\
& = \ringvol( (X_0+T)\cap \ringStrip C)  -\ringvol( ((X_0+T)\cap \ringStrip C) \cap C^+ ) \\
& = \ringvol(((X_0+T)\cap \ringStrip C)\cap C^-) \\
&=\frac{1}{2}
\end{align*}
as we wanted to show (cf. Figure \ref{fig:vw}, right).

\section*{Acknowledgements}

The author would like to thank Professor Takayubi Hibi and Akiyoshi Tuchiya for the invitation and amiable organization of the summer workshop on lattice polytopes 2018, and also all other participants for making it a wonderful and very productive stay.
Further thanks go to \emph{The Birth of Modern Trends on Commutative Algebra and Convex Polytopes with Statistical and Computational Strategies}
(Grant-in-Aid for Scientific Research (S) 26220701) for funding the workshop.

%%%%%%%%%%%%%%%%%%%%%%%%%%%%%%%%%%%%%%%%%%%%%%%%%%%%%%%%%%%%%%%%%%%%%%%%%%%%%%%%%
%%%%%%%%%%%%%%%%%%%%%%%%%%%%%%%%%%%%%%%%%%%%%%%%%%%%%%%%%%%%%%%%%%%%%%%%%%%%%%%%%

\bibliographystyle{ws-procs961x669}
%\bibliography{ws-pro-sample}

%\end{document}

%Non BiBTeX users can list down their references as:

\end{document}